\documentclass[letterpaper,conference]{IEEEtran}
\IEEEoverridecommandlockouts

\usepackage{hyperref}

\usepackage{amsfonts}
\usepackage{amsmath, amscd}
\usepackage{amsthm}
\usepackage{amssymb}
\usepackage{bbm, dsfont}
\usepackage[capitalize,noabbrev]{cleveref}
\usepackage[marginpar]{todo}
\usepackage{mathdots}
\usepackage{algorithm}
\usepackage{algorithmicx}
\usepackage{algpseudocode}
\usepackage{xcolor}
\usepackage{url}
\usepackage{optidef}
\usepackage{placeins}
\usepackage{graphicx}
\usepackage{multirow}
\usepackage{booktabs}
\usepackage{makecell}
\usepackage{notation}

\newtheorem{theorem}{Theorem}
\newtheorem{corollary}{Corollary}
\newtheorem{lemma}{Lemma}
\newtheorem{assumption}{Assumption}

\newcommand{\ti}{i}
\renewcommand{\vec}[1]{#1}

\title{Joint Problems in Learning Multiple Dynamical Systems\thanks{The final version of this paper was accepted to the 61st Allerton Conference on Communication, Control, and Computing, 2025, \url{https://allerton.csl.illinois.edu/}. This work was co-funded by the European Union under Horizon Europe grant number 101084642 and by the UK Research and Innovation (UKRI) under the UK government’s Horizon Europe funding guarantee [grant number 101084642], as part of the CoDiet project.}}

    \makeatletter
    \newcommand{\linebreakand}{%
      \end{@IEEEauthorhalign}
      \hfill\mbox{}\par
      \mbox{}\hfill\begin{@IEEEauthorhalign}
    }
    \makeatother

\author{
\IEEEauthorblockN{Mengjia Niu}
\IEEEauthorblockA{\textit{Imperial College London} \\
\textit{South Kensington Campus}\\
London SW7 2AZ, UK \\
\href{mailto:m.niu21@imperial.ac.uk}{m.niu21@imperial.ac.uk} }
\and
\IEEEauthorblockN{Xiaoyu He}
\IEEEauthorblockA{\textit{Czech Technical University in Prague} \\
\textit{Karlovo náměstí 13}\\
Prague, 12135, Czechia \\
\href{mailto:hexiaoyu@cvut.cz}{hexiaoyu@cvut.cz} }
\and
\IEEEauthorblockN{Petr Ryšavý}
\IEEEauthorblockA{\textit{Czech Technical University in Prague} \\
\textit{Karlovo náměstí 13}\\
Prague, 12135, Czechia \\
\href{mailto:petr.rysavy@fel.cvut.cz}{petr.rysavy@fel.cvut.cz} }
\linebreakand 
\IEEEauthorblockN{Quan Zhou}
\IEEEauthorblockA{\textit{Technion-Israel Institute of Technology} \\
\textit{Technion City}\\
Haifa, 3200003, Israel \\
\href{mailto:quan.zhou@campus.technion.ac.il}{quan.zhou@campus.technion.ac.il} }
\and
\IEEEauthorblockN{Jakub Mareček}
\IEEEauthorblockA{\textit{Czech Technical University in Prague} \\
\textit{Karlovo náměstí 13}\\
Prague, 12135, Czechia \\
\href{mailto:jakub.marecek@fel.cvut.cz}{jakub.marecek@fel.cvut.cz} }
}

\begin{document}
\maketitle

\begin{abstract}
Clustering of time series is a well-studied problem \cite{vlachos2002discovering,keogh2004indexing}, with applications ranging from quantitative, personalized models of metabolism obtained from metabolite concentrations to state discrimination in quantum information theory. We consider a variant, where given a set of trajectories and a number of parts, we jointly partition the set of trajectories and learn linear dynamical system (LDS) models for each part, so as to minimize the maximum error across all the models. We present globally convergent methods and EM heuristics, accompanied by promising computational results. The key highlight of this method is that it does not require a predefined hidden state dimension but instead provides an upper bound. Additionally, it offers guidance for determining regularization in the system identification. 
\end{abstract}

\section{Introduction}
The task of clustering similar time series based on their dynamic patterns has attracted significant attention due to its applications ranging from studying mobility patterns \cite{mokbel2023towards} to improving Apple Maps \cite{chen2016city}, 
through quantitative, personalized models of metabolism obtained from metabolite concentrations,
all the way to state discrimination problems in quantum information theory \cite{magesan2015machine}. 

We consider a variant, where given a set of trajectories and a number of parts, we jointly partition the set of trajectories and learn the autonomous discrete-time Linear Dynamical System (LDS) \cite{west2006bayesian} models:
\begin{equation}
\begin{aligned}
    \vec{\phi}_{t} &= \mathbf{G} \vec{\phi}_{t-1} + \textrm{normally distributed noise}, \\
    \vec{x}_{t} &= \mathbf{F'} \vec{\phi}_t + \textrm{normally distributed noise},
\end{aligned}
\label{eq:ldsmotivation}
\end{equation}
for each cluster, where $\vec{\phi}_t$ are the hidden states and $\vec{x}_t$ are the observations.
The cluster-specific LDSs may exhibit similar behaviors in terms of system matrices $\mathbf{F},\mathbf{G}$, or not. $\mathbf{F'}$ denotes the transpose of $\mathbf{F}$.
The observations convey information about the cluster-specific LDSs. 

The main contributions of this paper are the following.
\begin{itemize}
    \item We propose a novel problem in the clustering time-series considering Linear Dynamic System \cite{west2006bayesian} models for each cluster. The linearity assumption comes without a loss of generality as any non-linear system can be modeled as an LDS \cite{BEVANDA2021197}, in a sufficiently higher dimension. 
    \item We provide an abstract formulation within Non-Commutative Polynomial Optimization (NCPOP). NCPOP \cite{ncpop} is a framework for operator-valued optimization problems, and thus does not require the dimension of the hidden state to be known ahead of time, which had been known \cite{liu2015regularized} to be a major limitation of LDS-based methods. This paper is one of the first applications of NCPOP in machine learning. 
    \item As a complement to the NCPOP formulation, we provide an efficient Expectation-Maximization (EM) procedure \cite{emalg}.
    Through iterative measurement of prediction errors and systematic updates to the system matrix, we can effectively identify the per-cluster LDSs and the assignment of time series to clusters.
\end{itemize}


\section{Background}
\label{sec:background}

This section provides an overview and necessary definitions of the background, problems, and algorithms. A comprehensive table of notation is provided in the Supplementary Material.

\subsection{Linear Dynamic Systems (LDS) and System Identification}
\label{sec:systemid}

There is a long research procedure in system identification \cite{ljung2010perspectives} and related approaches in Bayesian statistics \cite{west2006bayesian}. 
Let $n$ be the hidden state dimension and $m$ be the observational dimension. A linear dynamic system $\mathbf{L}$ is defined as a quadruple $ (\mathbf{G}, \mathbf{F}, \mathbf{\Sigma}_H, \mathbf{\Sigma}_O)$, where $\mathbf{G}$ and $\mathbf{F}$ are \emph{system matrices} of dimension $n \times n$ and $n \times m$, respectively. Besides, $\mathbf{\Sigma}_H \in \mathbb{R}^{n \times n}$ and $\mathbf{\Sigma}_O \in \mathbb{R}^{m \times m}$ are covariance matrices \cite{west2006bayesian}. A single realization of the LDS or a \emph{trajectory} of length $T$ can be denoted by $\mathbf{X} = \{x_1, x_2, \ldots, x_T\} \in \mathbb{R}^{m \times T}$. And based on \emph{initial conditions} $\vec{\phi}_0$, and realization of noises $\vec{\upsilon}_t$ and $\vec{\omega}_t$, it is defined as 
\begin{align}
\phi_t &= \mat{G}\phi_{t-1} + \omega_t, \label{eq:lds1} \\
x_t &= \mat{F}' \phi_t  + \upsilon_t,
\label{eq:lds2}
\end{align}
where $\phi_t \in \mathbb{R}^{ n }$ is the vector autoregressive processes with hidden components. $\left\{\omega_t,\upsilon_t\right\}_{t \in \{1, 2, \ldots, T\}}$ are normally distributed process and observation noises with zero mean and covariance of $\mat{\Sigma}_H$ and $\mat{\Sigma}_O$ respectively, i.e., $\omega_{t}\sim N(0, \mat{\Sigma}_H)$ and  $\upsilon_{t}\sim N(0, \mat{\Sigma}_O)$. The transpose of $\mat{
F}$ is denoted as $\mat{F}'$. Vector $x_t \in \mathbb{R}^m$ serves as an observed output of the system. Recently, Zhou and Marecek \cite{zhou2020proper} proposed to find the global optimum of the objective function subject to the feasibility constraints arising from \eqref{eq:lds1} and \eqref{eq:lds2}:
\begin{equation}
\min_{f_t,\phi_t,\mat{G},\mat{F},\omega_t,\upsilon_t} \sum_{t \in \{1,2, \ldots, T\}} \|X_t - f_t \|_2^2 + \|\omega_t\|_2^2 + \|\upsilon_t\|_2^2,
\end{equation}
for a $L_2$-norm $\|\cdot\|_2$. In the joint problem we are given $N$ trajectories $\mat{X} \in \mathbb{R}^{m \times T}$. A natural problem to solve is to find the parameters of the LDS that generated the trajectories. In other words, we are interested in finding the optimal objective values, as well as system matrices $\mat{G}, \mat{F}$, and the noise vectors $\upsilon_t$, $\omega_t$ that belong to each LDS.

One should like to remark that learning the LDS is an NP-Hard problem. This is easy to see when one realizes \cite{roweis1999unifying} that Gaussian mixture models (GMM), autoregressive models, and hidden Markov models are all special cases of LDS, 
and whose learning is all NP-Hard, even in very special cases such as spherical Gaussians \cite{tosh2017maximum} in a GMM.
Furthermore, there are also inapproximability results \cite{tosh2017maximum} suggesting that there exists an approximation ratio,
at which no polynomial-time algorithm is possible unless P = NP.

\subsection{Clustering with LDS Assumptions}

The problem of clustering of time series is relevant in many fields, including \cite{reviewmulti,WARRENLIAO20051857} applications in Bioinformatics, Multimedia, Robotics, Climate, and Finance. There are a variety of existing methods, including those based on (Fast) Discrete Fourier Transforms (FFTs), Wavelet Transforms, Discrete Cosine Transformations, Singular Value Decomposition, Levenshtein distance, and Dynamic Time Warping (DTW).
We compare our method with the FFT- and DTW-based methods.

Many methods combine the ideas of system identification and clustering, sometimes providing tools for clustering time series generated by LDSs, similar to our paper. With three related papers at ICML 2023, this could be seen as a hot topic. 
We stress that neither of the papers has formulated the problem as either a mixed-integer program or an NCPOP,
or attempted to solve the joint problem without decomposition into multiple steps, 
which necessarily restricts both the quality of the solutions one can obtain in practice, 
as well as the strength of the guarantees that can be obtained in theory. 

To our knowledge, the first mention of clustering with LDS assumptions is in the paper of \cite{li2011time}, who introduced ComplexFit, a novel EM algorithm to learn the parameters for complex-valued LDSs and utilized it in clustering.
In \cite{liu2015regularized}, regularization has been used in learning linear dynamical systems for multivariate time series analysis.
In a little-known but excellent paper at AISTATS 2021, Hsu et al. \cite{hsu2020linear} consider clustering with LDS assumptions, 
but argue for clustering on the eigenspectrum of the state-transition matrix ($\mathbf{G}$ in our notation), which can be identified for unknown linear dynamics without full system identification.
The main technical contribution is bidirectional perturbation bounds to prove
that two LDSs have similar eigenvalues if and only if
their output time series have similar parameters within Autoregressive-Moving-Average (ARMA) models. Standard consistent estimators for
the autoregressive model parameters of the ARMA models are then used to estimate the
LDS eigenvalues, allowing for linear-time algorithms. 
We stress that the eigenvalues may not be interpretable as features; one has to provide the dimension of the hidden state as an input.  

At ICML 2023, Bakshi et al.\ \cite{pmlr-v202-bakshi23a} presented an algorithm to learn a mixture of linear dynamic systems. Each trajectory is generated so that an LDS is selected on the basis of the weights of the mixture, and then a trajectory is drawn from the LDS. Their approach is, unlike ours, based on moments and the Ho-Kalman algorithm and tensor decomposition, which is generalized to work with the mixture.
Empirically, \cite{pmlr-v202-bakshi23a} outperforms the previous work of Chen et al.\ \cite{chen} of the previous year, which works in fully observable settings. In the first step, the latter algorithm \cite{pmlr-v202-bakshi23a} finds subspaces that separate the trajectories. In the second step, a similarity matrix is calculated, which is then used in clustering and consequently can be used to estimate the model parameters. 
The paper \cite{pmlr-v202-bakshi23a} also discusses the possibility of classification of new trajectories and provides guarantees on the error of the final clustering.



There are also first applications of the joint problems in the domain-specific literature. Similarly to the previous paper, a fully observable setting of vector autoregressive models is considered in \cite{varmodels}, with applications in Psychology, namely on depression-related symptoms of young women. Similarly to our method, the least-squares objective is minimized to provide clustering in a manner similar to the EM-heuristic. See also \cite{Ernst2023Mixture} for further applications in Psychology. 
One can easily envision a number of further uses across Psychology and Neuroscience, especially when the use of mixed-integer programming solvers simplifies the time-consuming implementation of EM algorithms. 

\subsection{Mixed-Integer Programming (MIP)}
To search for global optima, we developed relaxations to bound the optimal objective values in non-convex Mixed-Integer Nonlinear Programs (MINLPs) \cite{belotti2013mixed}. Our study is based on the Mixed-Integer Nonlinear Programs of the form:
\begin{align}
\label{minlp}
\tag{MINLP}  \min f(x,y,z)\\
\textrm{s.t. } g_i(x,y,z)\leq 0, &\quad \forall i \in I,\\
h_j(x)\leq 0, &\quad \text{if } z_j= 1, \forall j \in J,\\
              &\quad x\in \mathbb{R}^{n},y\in \mathbb{Z}^{m}
\notag
\end{align}
where functions $f$, $g_i$ and $h_j$ are assumed to be continuous and twice differentiable.
Such problems are non-convex, both in terms of featuring integer variables and in terms of the functions $f, g_i$. 

While MINLP problems may seem too general a model for our joint problem, notice that the NP-hardness and inapproximability results discussed in Section~\ref{sec:systemid} suggest that this may be the appropriate framework. For bounded variables, standard branch-and-bound-and-cut algorithms 
run in finite time. Both in theory -- albeit under restrictive assumptions, such as in \cite{dey2021branch} -- and in practice, the expected runtime
is often polynomial. 
In the formulation of the next section, $f, g_i$ are trilinear, and various monomial envelopes have been considered and implemented in global optimization solvers such as BARON \cite{sahinidis1996baron}, SCIP \cite{bestuzheva2023global}, and Gurobi. In our approach, we consider a mixed-integer programming formulation of a piecewise polyhedral relaxation of a multilinear term using its convex-hull representation.

\subsection{Non-Commutative Polynomial Optimization (NCPOP)}

To extend the search for global optima from a fixed finite-dimensional state to an operator in an unknown dimension, we formulate the problem as a non-commutative polynomial optimization problem (NCPOP), cf. \cite{pironio2010convergent, burgdorf2016optimization}.
In contrast to traditional scalar-valued, vector-valued, or matrix-valued optimization techniques, the variables considered in NCPOP are operators, whose dimensions are unknown \textit{a priori}.


Let $X=(X_1,\ldots,X_n)$ be a tuple of bounded operators on a Hilbert space $\mathcal{H}$. Let $[X,X^{\dag}]$ denote these $2n$ operators, with the $\dag$-algebra being conjugate transpose. 
Let monomials $\omega,\mu$ be products of powers of variables from $[X,X^{\dag}]$.
The degree of a monomial, denoted by $|\omega|$, refers to the sum of the exponents of all operators in the monomial $\omega$, e.g., $|X^3_nX_n^{\dag}|=4$.
Let $p$ and $q_i$, $i=1,\dots,m$ be polynomials in these $2n$ variables. 
Let $\deg(p)$ denote the polynomial degree of $p$. 
In the following, we will view these
$2n$ variables as the new tuple $X$.
Using the set of monomials generated from the tuple $X$, polynomials $p$ and $q_i$, $i=1,\dots,m$ can be rewritten as linear combinations of monomials:
\begin{equation*}
    p(X)=\sum_{|\omega|\leq \deg(p)} p_{\omega} \omega,\quad
    q_i(X) = \sum_{|\mu|\leq \deg(q_i)} q_{i,\mu} \mu,
\end{equation*}
for $i = 1,\ldots,m$, and $p_{\omega},q_{i,\mu}$, are coefficients of these polynomials. For instance, $p(X)=X_1^3X_n^{\dag}+5X_n=\omega_1+5\omega_2$, where $\omega_1=X_1^3X_n^{\dag}$ and $\omega_2=X_n$.

Let $\langle\cdot,\cdot\rangle$ denotes inner product.
Suppose there is a normalized vector $\psi$, i.e., $\langle\psi,\psi\rangle=1$, also defined on the Hilbert space $\mathcal{H}$.
Let $p(X),q_i(X)$ be the Hermitian operators, i.e., $p^{\dag}(X)=p(X)$.
The formulation considered in NCPOP reads
\begin{mini}
{(\mathcal{H},X,\psi)} {\langle\psi,p(X)\psi\rangle}{}{}
\addConstraint{q_i(X)}{\succcurlyeq 0, }{\;i=1,\ldots,m}
\addConstraint{\langle\psi,\psi\rangle}{= 1,}{}
\end{mini}
where the constraint $q_i(X)\succcurlyeq 0$ denotes that the variable $q_i(X)$ is positive semidefinite. 



Under the Archimedean assumption, such that the tuple of operators $X$ are bounded, one can utilize the Sums of Squares theorem of \cite{helton2002positive} and \cite{mccullough2001factorization} to derive semidefinite programming (SDP) relaxations of the Navascules-Pironio-Acin (NPA) hierarchy \cite{Navascues,Pironio2010}.
There are also variants \cite{wang2019tssos,wang2020chordal,wang2020sparsejsr} that exploit various forms of sparsity. 

\section{Problem Formulation}
\label{sec:problem}
Suppose that we have $N$ trajectories denoted by $\mathbf{X}^\ti\in \mathbb{R}^{m \times T}$ for $i \in \{1, 2, \ldots, N\}$. 
We assume that those trajectories were from $K(=2)$ clusters, $C_0$ and $C_1$. The trajectory of the cluster $C_0$ (resp. $C_1$) was generated by a LDS $\mathbf{L_0}$ (resp. $\mathbf{L_1}$).
We aim to jointly partition the given set of trajectories into $K(=2)$ clusters and recover the parameters of the LDSs systems of both clusters $C_0,C_1$.
To solve these joint optimization problems, we can introduce an \emph{indicator function} to determine how the $N$ trajectories are assigned to two clusters:
\begin{equation}
\label{ln}
l_{\ti}=\begin{cases}
0, & \text{if}\ \ti \in C_0,\\
1, & \text{if}\ \ti \in C_1,
\end{cases}
\end{equation}
for $\ti \in \{1, 2, \ldots, N\}$. 

\subsection{Least-Squares Objective Function} \label{sec:optimization}

We define the optimization problem with a least-squares objective that minimizes the difference of measurement estimates $\vec{f}^0_t \in \mathbb{R}^{n\times T}$, $\vec{f}^1_t \in \mathbb{R}^{n\times T}$ and the corresponding measurements.
Other variables include noise vectors that come with the estimates; indicator function $l_i$ that assigns the trajectories to the clusters, and parameters of systems $\mathbf{L_0}$ and $\mathbf{L_1}$ (also known as system matrices). The objective function reads:
\begin{multline}
    \min_{
      \substack{ l_i\in \{0, 1\}
    }}
        \sum_{\ti=1}^N \sum_{t=1}^{T} 
        \|\vec{X}^\ti_t - \vec{f}^{l_i}_t\|_2^2 + 
       \sum_{c \in \{0, 1\}}\sum_{t=1}^{T} \left[  \| \vec{\upsilon_{t}}^{c} \|_2^2 + \| \vec{\omega_{t}}^{c} \|_2^2 \right]\\
\omega_{t}^{c}\sim N(0, \mat{\Sigma}_H^{c}), 
\upsilon_{t}^{c}\sim N(0, \mat{\Sigma}_O^{c}),
    \label{eq:optcriterion}
\end{multline}
where $\vec{\upsilon}^{c}, \vec{\omega}^{c}$ are defined above 
and $\| \cdot \|_2$ denotes the $L_2$ norm. Note that the indicator index $l_i$ in the superscript in \eqref{eq:optcriterion} can be replaced by multiplication with the indicator function, i.e., 
\begin{align}
    \|\vec{X}^\ti_t - \vec{f}^{l_i}_t\|_2^2 = \|\vec{X}^\ti_t - \vec{f}^{0}_t\|_2^2 \cdot (1-l_i) + \|\vec{X}^i_t - \vec{f}^{1}_t\|_2^2 \cdot l_i.
    \label{eq:indicatorbrevity}
\end{align}
In the first part of the optimization criterion \eqref{eq:optcriterion}, we have a sum of squares of the difference between trajectory estimate $\vec{f}_t^{l_\ti}$ and observations of the trajectories assigned to cluster $C_{c}$ for $c=0,1$. The second part of the optimization criterion \eqref{eq:indicatorbrevity} 
provides a form of regularization with parameter=$270$.

\subsection{Feasible Set in State Space}

The feasible set given by constraints is as follows:
\begin{align}
    \vec{\phi}^{c}_t &= \mat{G}_{c} \vec{\phi}^{c}_{t-1} + \vec{\omega}^{c}_t,         && \forall t=2,\dots,T, && \forall c=0,1, \label{eq:constr1}\\
    \vec{f}^{c}_t &= \mat{F}^{'}_{c} \vec{\phi}^{c}_t + \vec{\upsilon}^{c}_t, && \forall t=1,\dots,T,&& \forall c=0,1, \label{eq:constr2} \\
    l_\ti^2 & = l_\ti && \forall \ti=1,\dots,N. \label{eq:constr3}
\end{align}
The first two equations in the set of constraints, \eqref{eq:constr1} and \eqref{eq:constr2}, encode that the system is an LDS with system matrices $\mat{F}$ and $\mat{G}$. For brevity of the notation, both equations are indexed by the cluster index $c$, which can be rewritten as twice as many equations, one with $\vec{f}^0_t, \mathbf{F}_0, \mathbf{G}_0, \vec{\upsilon}^0_t, \vec{\omega}^0_t$, and $\vec{\phi}_t^0$, the second one with $\vec{f}^1_t, \mathbf{F}_1, \mathbf{G}_1, \vec{\upsilon}^1_t, \vec{\omega}^1_t$, and $\vec{\phi}_t^1$. The third equation \eqref{eq:constr3} encodes that the indicator function is $0$ or $1$ for each trajectory.

A weighted combination of the redundant constraints in the spirit of Gomory. 
While these strengthen the relaxations, the higher-degree polynomials involved come at a considerable cost. 
Still, even when the dimension of the hidden state \vec{n} is unknown, we can solve the corresponding operator-valued problem:

\begin{theorem}
There exists a series of convex relaxations, whose optima asymptotically converge to the true global optimizer of the problem Equation~\eqref{eq:optcriterion} subject to (\ref{eq:constr1}--\ref{eq:constr3}).
\end{theorem}
\begin{proof}
Let $n,m$ be positive integers,
$x\in\mathbb{R}^n$ be a tuple of real-valued variables, and $p,q_i$, $i=1,\ldots,m$ be some polynomials in the variable $x$.
Polynomial optimization consider $\min_{x\in\mathbb{R}^n} p(x)$ subject to $q_i(x)\geq 0$ for $i=1,\ldots,m$.
Under the Archimedeam assumption, such feasible region is a compact semi-algebraic set. 
Note that the formulation is equivalent to finding the maximum number $\alpha$ that makes the polynomial $p(x)-\alpha$ nonnegative on the compact semi-algebraic set defined by $q_i(x)$, $i=1,\ldots,m$.
Then, according to the Putinar's certificate of positive polynomials (i.e., Putinar's positivstellensatz), if a polynomial $p$ is strictly positive on a compact semi-algebraic set, there exists a sequence of sum-of-square polynomials $g_i$, $i=0,\ldots,m$ such that $p=g_0+\sum_{i=1}^{m} q_i g_i$, where verifying sum-of-square polynomials is by solving SDP problems.
Considering this, Lasserre's hierarchy of SDP relaxations provides global convergence for Polynomial optimisation \cite{lasserre2001global,lasserre2009moments}, following Putinar's positivstellensatz and Curto-Fialkow’s theorem.

NCPOP is the extension of polynomial optimization to consider the variables $X=\{X_1,\dots, X_n\}$ which are not simply real numbers but non-commutative variables, for which, in general, $X_i X_j\neq X_j X_i$. 
The polynomials e.g., $p(x)$, are defined by substituting the variables $x$ by the tuple of operators $X$ in the expression $p(X)$.
The global convergence was provided in Navascu\'es-Pironio-Ac\'in (NPA) hierarchy of SDP relaxations \cite{pironio2010convergent,navascues2012sdp}, following Helton and McCullough's certificate of non-commutative positive polynomials \cite{helton2004positivstellensatz}. 
To lower the computational burden of NCPOP, the sparsity-exploiting variants were provided \cite{wang2019tssos,klep2021sparse}.
\end{proof}

Despite the existence of the relaxations, we can show that the soft-clustering version of the problem is NP-hard, as the problem can be transformed to finding a clustering of a mixture of Gaussian distributions, a related and well-studied problem known to be NP-hard even for spherical clusters \cite{gaussiannphard}.
\begin{theorem}
    Finding a soft clustering of a mixture of LDS trajectories with a log-likelihood within an additive factor of the optimal log-likelihood is NP-hard even when $k = 2$.
\end{theorem}

\subsection{Variants and Guarantees}
\label{sec:variants}

There are several variants of the formulation above. Notably, one can:
\begin{itemize}
\item consider a fixed, finite dimension of the hidden state $\phi$ to be known and to solve a finite-dimensional \eqref{minlp}.
\item consider side constraints on the system matrices $\mat{F}, \mat{G}$, as in Ahmadi and El Khadir \cite{ahmadi2020learning} -- or not. 
At least requiring the norm of $\mat{G}$ to be $1$ is without loss of generality. 
\item bound the magnitude of the process noise $\vec{\omega}^{c}$ and observation noise $\vec{\upsilon}^{c}$, or other shape constraints thereupon.
\item bound the cardinality of the clusters -- or not.
\end{itemize}

Throughout, one obtains asymptotic guarantees on the convergence of the NCPOP, 
or guarantees of finite convergence in the case of the MINLP. 

\section{EM Heuristic}
\label{sec:em}

In addition to tackling the optimization problem in Section \ref{sec:optimization}, we provide an alternative solver using the Expectation-Maximization (EM) heuristic \cite{emalg}. The main idea of the algorithm is to avoid the direct optimization of the criterion in \eqref{eq:optcriterion}. Instead, the indicator function is treated separately. In the expectation step, the parameters of the LDSs are fixed, and the assignment of the trajectories to the clusters (i.e., the indicator function) is calculated. Then, in the maximization step, the criterion is optimized, and the LDS parameters are calculated with the indicator function fixed.
See Algorithm \ref{alg:emheuristic} for the pseudocode. First, the algorithm randomly partitions the trajectories into the clusters. Then, for each cluster, the parameters of the LDSs are found, and with the parameters known, the optimization criterion is recalculated, and each trajectory is put to the cluster, which lowers the error in \eqref{eq:optcriterion}.

The advantage is that the problem of finding the parameters and an optimal trajectory for a set of trajectories is easier than clustering the trajectories. 
Another advantage of the EM heuristic is that it can be easily generalized to an arbitrary number of clusters $K$, generally for any $K > 1$. 

In the supplementary materials, we prove the following theorem that shows that the problem of clustering a mixture of LDSs is no more difficult than clustering a mixture of Gaussian distributions as below.
\begin{theorem}
    \label{thm:reduction}
    There exists a polynomial reduction that reduces the problem of clustering a mixture of autonomous LDSs with hidden states to the clustering of a mixture of Gaussian distributions.
\end{theorem}
The theorem justifies the usage of the EM-algorithm. Unfortunately, applying the previous theorem directly to the joint problem is computationally inefficient, as a quadratic number of parameters needs to be estimated. The advantage is that we can exploit the theoretical guarantees for the mixture of Gaussians problem - for example, the local convergence to a global optimum \cite{zhao2020}, existence of arbitrarily bad local minima \cite{NIPS2016_3875115b}, and also a linear bound on the number of samples in the case of spherical clusters \cite{pmlr-v125-kwon20a}. In the case when there are only two clusters, the EM-algorithm-based estimates are guaranteed to converge to one of three cases \cite{NIPS2016_792c7b5a}.  See the Supplementary materials for more details.

It would be of considerable interest to analyse the behaviour of the EM heuristic in our setting. 
Indeed, for many problems, such as the parameter estimation of Gaussian mixture models 
\cite{dwivedi2020sharp,weinberger2022algorithm}, the properties of EM approaches are well understood \cite{dwivedi2020sharp,ho2022weak,weinberger2022algorithm,zhang2022distributed,ho2022convergence}.
The joint problems are very similar to the clustering of mixture of Gaussian distributions over the system matrices as we have seen in Theorem \ref{thm:reduction} - in our setting of autonomous LDSs with hidden state, one can treat all observations of a trajectory as individual features in a high-dimensional space and the resulting vector will follow a normal distribution with additional constraints applied on its parameters.

As EM-algorithm applied to the mixture of Gaussians is, in many scenarios, computationally inefficient, we propose to use heuristic \ref{alg:emheuristic}, which can be seen as a parallel to the Lloyd's algorithm \cite{lloydkmeans} for $k$-means problem. See the Supplementary materials for a formal connection to the $k$-means problem.

\begin{algorithm}[t!]
\begin{algorithmic}
\Function{EM-clustering}{$N$ trajectories $\mat{X}^\ti\in \mathbb{R}^{m \times T}$, $K$}
  \State\Comment{Generate a random partitioning into two clusters.}
  \State $\vec{l}_i \gets \Call{RandomInt}{\{0,1,\ldots, K - 1\}}$
  \State
  \State\Comment{Iterate until convergence.}
  \While{$l_i$ changes for any $i \in \{1,2, \ldots, N\}$}
     \State\Comment{For each cluster find cluster parameters}
     \For{$c \in \{0,1,\ldots, K\}$}
        \State Find the cluster $C_{c}$ parameters by solving 
        \begin{multline*}
            \min_{
              \substack{
               \vec{f}^c_t, \mat{F}_c, \mat{G}_c, \vec{\upsilon}^c_t, \vec{\omega}^c_t, \vec{\phi}_c^0
            }}
              \left[ \sum_{\ti=1}^N \sum_{t=1}^{T} 
                \mathbbm{1}[l_i = c] \cdot
                \|\vec{X}^\ti_t - \vec{f}^{c}_t\|_2^2 \right] \\ +
                \| \vec{\upsilon}^{c} \|_2^2 + \| \vec{\omega
        }^{c} \|_2^2
        \end{multline*}
     \EndFor
    \State
     \State\Comment{Reassign the trajectories to the clusters.}
     \For{$\ti \in \{1, 2, \ldots N\}$}
       \State $l_\ti \gets \displaystyle\mathop{\mathrm{arg\, min}}_{c \in \{0,1,\ldots, K-1\}} \sum_{t=1}^T \|X^\ti_t - f^c_t\|_2^2 $
     \EndFor
  \EndWhile
\EndFunction
\end{algorithmic}
\caption{The EM heuristic.}
\label{alg:emheuristic}
\end{algorithm}

\section{Experiments}
\label{sec:experiments}

In this section, we present a comprehensive set of experiments to evaluate the effectiveness of the proposed \eqref{minlp} and EM heuristic,
without considering any side constraints and without any shape constraints.  
Our experiments are conducted on Google Colab with two Intel(R) Xeon(R) (2.20GHz) CPUs and Ubuntu 22.04.
The source code is included in the Supplementary Material.


\subsection{Methods and Solvers}
\label{sec:general settings}
\paragraph{MIP-IF}
Our formulation for MINLP with an indicator function (MIP-IF) uses equations \eqref{eq:indicatorbrevity} subject to (\ref{eq:constr1}--\ref{eq:constr3}). We specify the dimension of system matrix $\mat{G}$, i.e., $n$ as the hyperparameter. For every data set, we perform $50$ experiments, in which different random seeds are employed to initialize the indicator function in each iteration. These MIP instances are solved via Bonmin\footnote{https://www.coin-or.org/Bonmin/} \cite{bonami2008algorithmic} and Gurobi \footnote{https://www.gurobi.com/} \cite{gurobi}.

\paragraph{EM Heuristic}
Iterated EM Heuristic clustering is presented in Algorithm~\ref{alg:emheuristic}. As in MIP-IF, the upper bound of dimension, $n$, of the system matrix is required. In each iteration, the cluster partition is randomly initialized and we conduct $50$ trials with various random seeds for every dataset. 

\begin{figure*}[t]
\centering
\includegraphics[width=0.98\textwidth]{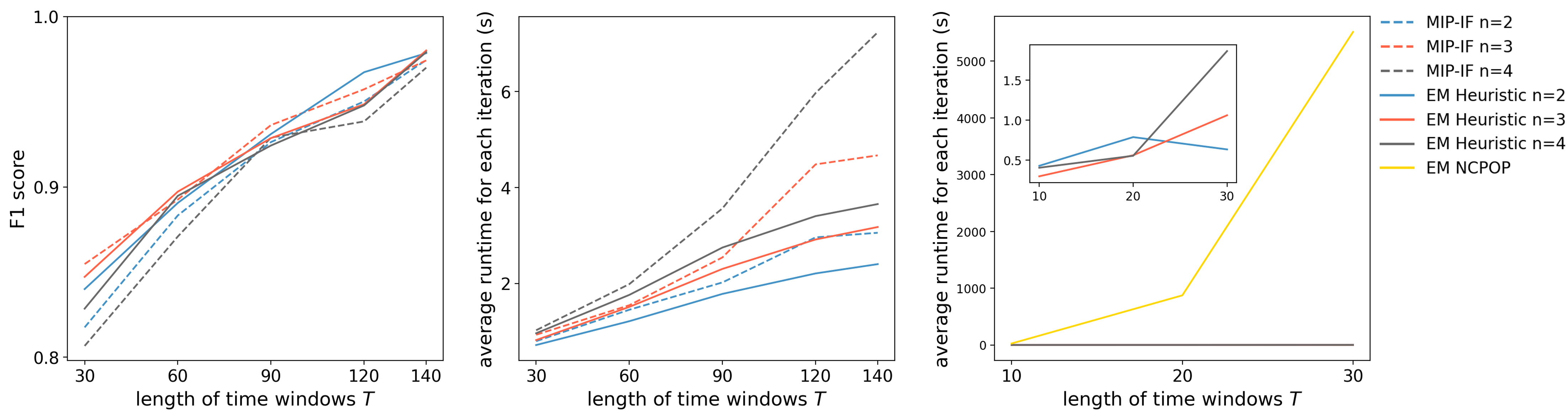}
\caption{\textbf{left} and \textbf{center}: $F_1$ scores and runtimes of MIP-IF and EM Heuristic with various time window $T$ chosen from $\{30,60,90,120,140\}$ respectively. The $F_1$ score improves as the time window increases. \textbf{right}: Runtimes of EM Heuristic using NCPOP compared with the method requiring specific dimension $n$ of hidden state. $n$ is chosen from $\{2,3,4\}$.}
\label{fig:ecg-mip-em-ncpop}
\end{figure*}

\textbf{The discussion of runtime}
As a subroutine of the EM Heuristic, the LDS identification in equation \eqref{eq:optcriterion} is solved via Bonmin, with runtime presented in the center subplot of Figure~\ref{fig:ecg-mip-em-ncpop}.
When the dimensions of system matrices are not assumed, the LDS identification becomes an NCPOP, and the runtime increases exponentially if NPA hierarchy \cite{Navascues,Pironio2010} is used to find the global optimal solutions, but stays relatively modest if sparity exploit variants  \cite{wang2019tssos,wang2020chordal,wang2020sparsejsr} are used.


\paragraph{Baselines}
We consider the following traditional time series clustering methods as baselines:
\begin{itemize}
    \item \textbf{Dynamic Time Warping (DTW)} is used for measuring similarity between given time series \cite{sakoe1978dynamic,gunopulos2001time}. We utilize K-means on DTW distance with tslearn \cite{JMLR:v21:20-091} package. 
    \item \textbf{Fast Fourier Transform (FFT)} is utilized to obtain the Fourier coefficients and the distance between time series is evaluated as the L2-norm of the difference in their respective Fourier coefficients. Subsequently, K-means is employed to cluster time series using this distance.
\end{itemize}

\begin{figure}[t]
\centering
\includegraphics[width=0.4\textwidth]{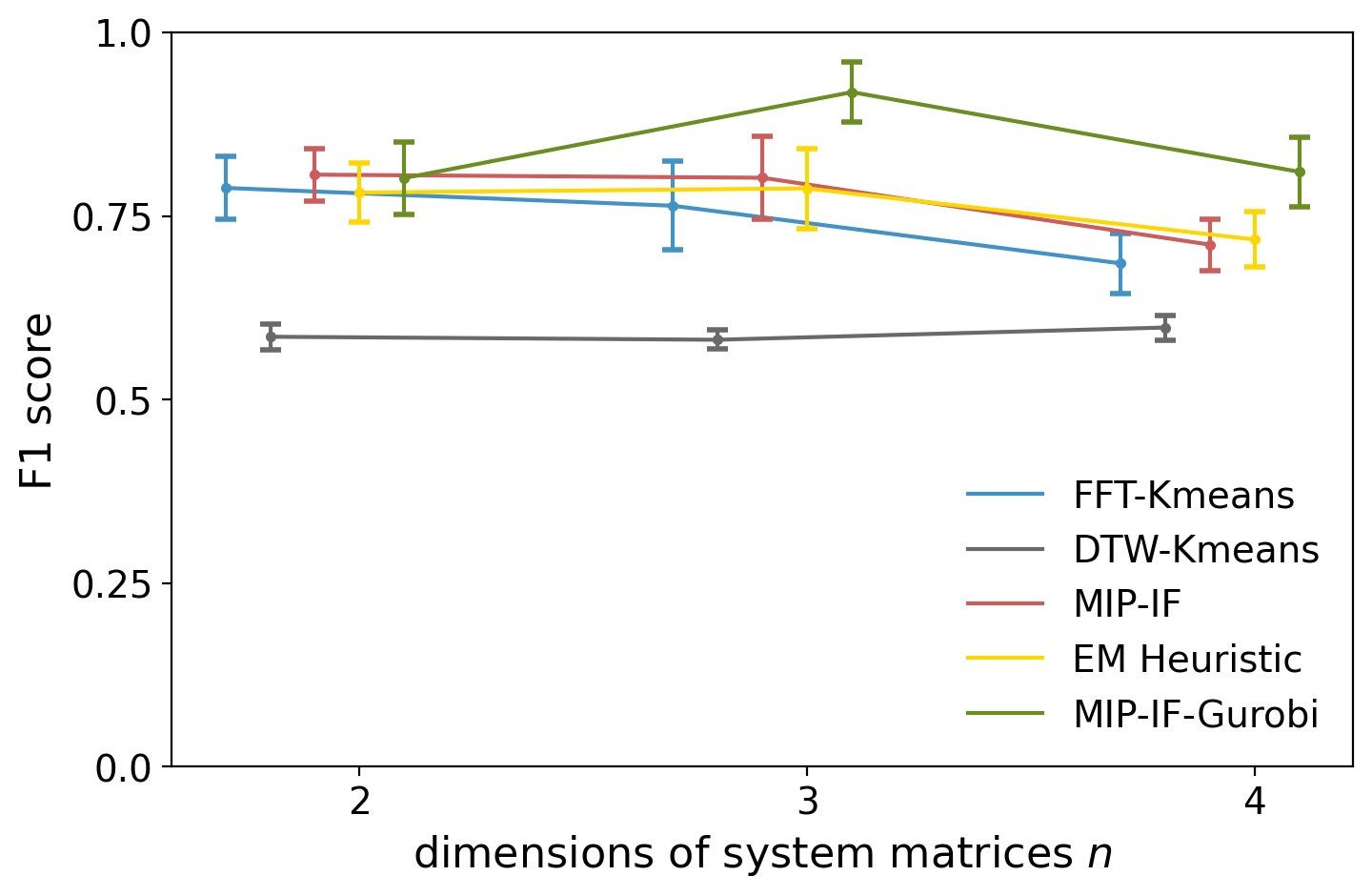}
\caption{$F_1$ scores of proposed methods compared to baselines. The results are based on the data generated by LDSs with hidden state dimensions $n \in \{2,3,4\}$ respectively. The vertical error bars show the $95\%$ confidence intervals of the $50$ trials. A higher $F_1$ score indicates better clustering performance.}
\label{fig:f1-N}
\end{figure}

\begin{figure}[t]
\centering
\includegraphics[width=0.4\textwidth]{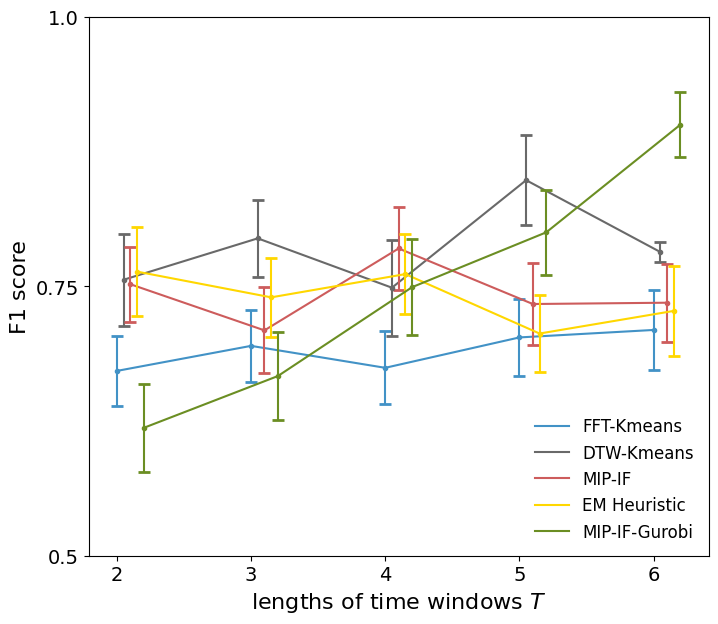}
\caption{$F_1$ scores of proposed methods compared to baselines. The results are based on data sampled from ECG5000 with various time window $T$ chosen from $\{30,60,90,120,140\}$. A higher $F_1$ score indicates better clustering performance.}
\label{fig:f1-ECG}
\end{figure}


\subsection{Experiments on Synthetic Data}

\paragraph{Data Generation}
\label{sec:synthdatasets}
We generate LDSs by specifying the quadruple $ (\mathbf{G}, \mathbf{F}, \mathbf{\Sigma}_H, \mathbf{\Sigma}_O)$ and the initial hidden state $\phi_0$. $\mathbf{G}$ and $\mathbf{F}$ are system matrices of dimensions $n \times n$ and $n \times m$, respectively. For each cluster, we derive $\frac{N}{2}$ trajectories including $T$ observations. The dimension $m$ of the observations is set to $2$ and $n$ is chosen from $\{2,3,4\}$. To make these LDSs close to the center of the respective cluster, we fix the system matrices while only changing the covariance matrices $\mathbf{\Sigma}_H, \mathbf{\Sigma}_O$ from $0.0004, 0.0016, 0.0036$, to $0.0064$ respectively. Note that $\mathbf{\Sigma}_H=0.0004$ here refers to $\mathbf{\Sigma}_H=0.0004\times\mathbf{I}_n$, where $\mathbf{I}_n$ is the $n\times n$ identity matrix. 
We consider all combinations of covariance matrices ($16$ trajectories) for each cluster. 

\paragraph{Results}
$F_1$ score is exploited to evaluate models' performance. See the Supplementary Material for further details. 

In our first simulation, we explore the effectiveness of the proposed methods in synthetic datasets. For each choice of $n \in \{2,3,4\}$, we run $50$ trials. In each trial, the indicator function $l_i$ for MIP-IF and the original clustering partition for EM Heuristic are randomly initialized.

Figure \ref{fig:f1-N} illustrates the $F_1$ score of our methods and baselines, with $95\%$ confidence intervals from $50$ trials. Different approaches are distinguished by colors.
Both solutions proposed yield superior cluster performance considering $n \in \{3,4\}$. When $n = 2$, our methods can achieve comparable performance to FFT-Kmeans. 
These experiments thus demonstrate the effectiveness of our approach. 

\subsection{Experiments on Real-world Data}
\label{sec:rwdatasetexperiments}
Next, we conduct experiments on real-world data. 

\paragraph{ECG data}
\label{sec:rwdatasets}
The test on electrocardiogram (ECG) data gives an inspiring application on guiding cardiologist's diagnosis and treatment \cite{Nezamabadi2023Unsupervised}. The ECG data
ECG5000 \cite{UCRArchive2018} is a common dataset for evaluating methods for ECG data, 
which has also been utilized by other papers \cite{hsu2020linear} on clustering with LDS assumptions.
The original data comes from Physionet \cite{baim1986survival,PhysioNet} and contains a $20$-hour-long ECG for congestive heart failure. 
After processing, ECG5000 includes 500 sequences, where there are $292$ normal samples and $208$ samples of four types of heart failure. Each sequence contains a whole period of heartbeat with $140$ time stamps. 

\paragraph{Results on ECG}

We randomly sample two clusters from normal sequences and one type of abnormal sequences respectively. As the entire period of time series data is not always available, we also extract subsequences with various lengths of time window $T$ chosen from $\{30,60,90,120,140\}$ to test the clustering performance. In Figure \ref{fig:f1-ECG}, with the assumption of the upper bound of hidden state dimension is $n=5$, we implement all methods for $50$ runs at each length of time window. 
Our methods exhibit competitive performance relative to FFT and DTW when $T=140$. When the time window decreases, the performance of the baselines significantly deteriorates, while our methods maintain a higher level of robustness.

In the left two subplots of Figure \ref{fig:ecg-mip-em-ncpop}, we further explore the performance of our methods at varying dimensions of the hidden state ($n\in\{2,3,4\}$), because the dimension of the hidden state $n$ of the ECG data is, indeed, unknown. 
When the length of the time window increases, both methods experience a slight improvement in clustering performance, but this performance remains relatively stable when the dimension $n$ changes.
The runtime is presented in the center subplot. 
Compared to MIP-IF, the EM Heuristic exhibits a modest growth in runtime as the length of the time window increases.

\begin{table}[ht]
\centering
\begin{tabular}{lccccc}
\toprule
\multicolumn{2}{l}{Methods} & EM n=2 & EM n=3 & EM n=4 & NCPOP \\
\midrule
\multirow{3}{*}{\makecell{F1 \\score}} & T=10 & 0.728 & 0.619 & 0.788 & \textbf{0.794} \\
  & T=20 & 0.805 & 0.897 & \textbf{0.927} & 0.699 \\
  & T=30 & 0.843 & 0.764 & 0.842 & \textbf{0.927} \\
\bottomrule
\end{tabular}
\caption{F1 score of EM Heuristic using NCPOP compared with the method requiring specific dimension $n$ of hidden state with window size $T$ in $\{10,20,30\}$. $n$ is chosen from $\{2,3,4\}$.}
\label{tab:ncpop}
\end{table}
Finally, when the dimension $n$ of the hidden state is not assumed, the subproblem of the EM heuristic becomes an NCPOP. 
For the implementation with such an assumption, we construct NCPOP using ncpol2sdpa 1.12.2\footnote{https://ncpol2sdpa.readthedocs.io/en/stable/} \cite{wittek2015algorithm}. Subsequently, the relaxation problem is solved by Mosek 10.1\footnote{https://www.mosek.com/} \cite{mosek2023mosek}. Noted that the execution time of NCPOP escalates rapidly as the trajectory length $T$ grows, we test NCPOP and compare its performance with the aforementioned EM Heuristic method with $T\in\{10,20,30\}$.
For comparison, we use pyomo\footnote{https://www.pyomo.org/} \cite{bynum2021pyomo} to construct the model and solve the problem with Bonmin\footnote{https://www.coin-or.org/Bonmin/} \cite{bonami2008algorithmic}, as above, with the dimension $n$ from $\{2,3,4\}$. 
The overall performance is illustrated in Table \ref{tab:ncpop}. NCPOP demonstrates the best performance in terms of the $F_1$ score when $T=10$ and $30$. 
However, as shown in the right subplots of Figure \ref{fig:ecg-mip-em-ncpop}, the runtime of NCPOP grows significantly as the length $T$ of the trajectory increases.

\section{Conclusions and Further Work}

We have studied problems in clustering time series, where given a set of trajectories and a number of parts, we jointly partition the set of trajectories and estimate a linear dynamical system (LDS) model for each part, so as to minimize the maximum error across all the models.
As discussed in Section \ref{sec:variants}, a number of variants of the joint problem remain to be investigated. The computational aspects of the operator-valued problem \cite{zhou2020proper} that consider the dimension of the hidden state to be unknown seem particularly interesting. 



\newpage
\bibliographystyle{plain}
\bibliography{refs,ncpop} 
\clearpage 

\newpage
\onecolumn
\setcounter{section}{0}

\begin{center}
\Large Supplementary Materials \\
\hrule
\end{center}

\section{Properties of Presented Method}
Considering the above contributions with previous work produced by \cite{li2011time,hsu2020linear,modi2022joint}, a pioneering method can be created that holds following properties:
\begin{itemize}
\item Mimic temporal dynamics in the data and produces the related, interpretable features, unlike methods relying on deep learning.
\item Handle multiple time lags within the time series without losing its periodicity, unlike methods based on Fourier coefficients \cite{kalpakis2001distance}. 
\item Not necessary to make any linear-algebraic assumptions, unlike subspace methods  \cite{agrawal1998automatic,domeniconi2004subspace} that assume the existence of a subspace, which might not necessarily be true. 
The assumption that the time series were generated by an LDS does not cause any loss of generality.
\item Avoid of making any assumptions on the dimension of the hidden state, unlike previous methods utilizing LDS \cite{li2011time,liu2015regularized,pmlr-v202-bakshi23a,modi2022joint}.
\item Improves upon previously proposed methods based on Fast Fourier Transform (FFT) \cite{kalpakis2001distance} and Dynamic Time Warping (DTW) \cite{sakoe1978dynamic,gunopulos2001time,Shah2016DTW} in practical experiments of Section \ref{sec:experiments}.
\end{itemize}

\section{Table of Notation Used in the Paper}

\begin{table}[h]
\caption{A Table of Notation.}
	\label{table:1}
\centering
\vskip 6pt
\begin{tabular}{| c | c | }
		\hline
		Symbol & Representation \\ \hline
            $i$ & Trajectory index \\ 
            $\mat{X}^i$ & Observation for a trajectory \\ 
            $t$ & Time index \\ 
            $T$ & Length of a trajectory \\ 
            $K$ & Number of clusters \\ 
            $N$ & Number of trajectories \\ 
            \hline
            $\phi_t$ & Hidden state processes \\ 
            $f_t$ & Optimal trajectory estimates \\ 
            $\omega_t$ & Hidden state noises \\ 
            $\upsilon_t$ & Observation error \\ 
		$n$ & Hidden state dimension \\ 
            $m$ & Observational dimension \\ 
            \hline
            $\mathbf{L}$ & LDS systems specified on $ (\mathbf{G}, \mathbf{F}, \mathbf{\Sigma}_H, \mathbf{\Sigma}_O)$   \\ 
		$\mat{G}$ & System matrix \\ 
            $\mat{F}$ & System matrix \\ 
            $\mat{\Sigma}
_H$ & Hidden state noises covariance matrix\\ 
            $\mat{\Sigma}_O$ & Observation error covariance matrix \\ 
            $l_i$ & Assignment of the trajectory $i$ \\ 
            $\vec{f}^{l_i}_t$ & Optimal trajectory estimates on ${l_i}$\\ 
            \hline
            $c$ & Cluser index \\ 
            $C_{c}$ & Set of trajectories assigned to c \\ 
            $\mathbf{L_c}$ & LDS systems generated by $C_{c} $ \\ 
		$\mat{G}_{c}$ & System matrix of $\mathbf{L_c}$\\ 
            $\mat{F}_{c}$ & System matrix of $\mathbf{L_c}$\\ 
            $\mat{\Sigma}_H^{c}$ & Hidden covariance matrix of $\mathbf{L_c}$\\ 
            $\mat{\Sigma}_O^{c}$ & Observation covariance matrix of $\mathbf{L_c}$\\ 
            $\phi_t^c$ & Hidden state processes of $\mathbf{L_c}$\\ 
            $\omega_{t}^{c}$ & Hidden state noises produced by $\mathbf{L_c}$\\
            $\upsilon_t^{c}$ & Observation error produced by $\mathbf{L_c}$\\ \hline
	\end{tabular}
\end{table}

\clearpage 

\section{Analysis of the EM-algorithm}

In this section, we will show that we can apply the EM-algorithm to the joint problem of clustering trajectories produced by multiple LDS while maintaining the properties of the EM-algorithm when applied to a mixture of Gaussian distributions. The overall idea is that we formalize the assumptions under which an autonomous linear dynamic system produces normally distributed observations at each time step. As the consecutive time steps are connected linearly, we will show that the resulting distribution will be Gaussian if we concatenate all time steps together in a single feature vector.

Through the text, assume that $\eye_m$ is an $\mathbb{R}^{m \times m}$ identity matrix.

\begin{assumption}
  \label{assumption:noisenormal}
  The hidden state noise $\hnoiset$ follows normal distribution $\normal(\mathbf{0}, \hcovariance)$. The observation noise $\onoiset$ follows $\normal(\mathbf{0}, \ocovariance)$.
\end{assumption}

\begin{assumption}
  \label{assumption:independence}
  Hidden state noise $\hnoiset$ and observation noise $\onoiset$ are both independent of the state/observation values and between their samples.
\end{assumption}

\begin{assumption}
  \label{assumption:initialstate}
  The hidden state $\hstate_0$ is normally distributed, i.e., $\hstate_0 \sim \normal( \hstatemuz, \hstatecovz )$.
\end{assumption}

\subsection{Preliminaries}

In the next section, we will need to use some well-known facts about the normal distribution and related consequences. We will formally state those preliminaries in this section.
\begin{lemma}[Linear transformation theorem of the multivariate normal distribution]
  \label{thm:lttmnd}
  Let
    $$x \sim \mathcal{N}(\mu, \Sigma).$$
  Then, any linear transformation of $x$ is also normally distributed
  $$ \mat{A}x + b \sim \mathcal{N}(\mat{A} \mu + b, \mat{A}\Sigma\mat{A}').$$
\end{lemma}

\begin{lemma}
  \label{thm:summnd}
  Let $x \sim \mathcal{N}(\mu_x, \Sigma_x)$ and $y \sim (\mu_y, \Sigma_y)$ be two independent, normally distributed multivariate normal distributions with $n$ dimensions. Then,
  $$ x + y \sim \mathcal{N}(\mu_x + \mu_y, \Sigma_x + \Sigma_y). $$
\end{lemma}
\begin{proof}
  Since $x$ and $y$ are independent, then
  $$ \mathcal{N}\left( \begin{bmatrix}
    \mu_x \\ \mu_y
  \end{bmatrix}, \begin{bmatrix}
    \Sigma_x & 0 \\ 0 & \Sigma_y
  \end{bmatrix} \right) $$
  is normally distributed. Using transformation matrix
  $$
    \mat{A} = \begin{bmatrix}
      \eye_n & \eye_n
    \end{bmatrix},
  $$
  where $\eye_n \in \mathbb{R}^{n \times n}$ is the identity matrix, the lemma is a direct result of Lemma \ref{thm:lttmnd}.
\end{proof}

\begin{lemma}
  \label{thm:gaussianonestep}
  Let $\vec{x} \sim \normal(\mu_x, \Sigma_x) \in \mathbb{R}^m$ be a normal distribution, and $\vec{y} \sim \normal(0, \Sigma_y) \in \mathbb{R}^n$ be an independent Gaussian noise. Then, concatenation of $\vec{x}$ and $\mathbf{A}\vec{x} + \vec{y}$ (where $\mathbf{A} \in \mathbb{R}^{n \times m}$) follows the normal distribution, i.e.,
  \begin{equation}
	\begin{pmatrix} \vec{x} \\ \mathbf{A}\vec{x} + \vec{y} \end{pmatrix}
	\sim
	\normal \left(
	  \begin{pmatrix}
	    \mu_x \\ \mathbf{A}\vec{x}
	  \end{pmatrix},
	  \begin{bmatrix}
	     \Sigma_x & \Sigma_x \mathbf{A}' \\
	     \mathbf{A} \Sigma_x & \mathbf{A} \Sigma_x \mathbf{A}' + \Sigma_y
      \end{bmatrix}	  	  
	\right).
  \end{equation}
\end{lemma}
\begin{proof}
As $\vec{x}$ and $\vec{y}$ are independent normal distributions, their concatenation is the following normal distribution
  \begin{equation}
    \begin{pmatrix} \vec{x} \\ \vec{y} \end{pmatrix}
	\sim
	\normal \left(
	  \begin{pmatrix}
	    \mu_x \\ 0
	  \end{pmatrix},
	  \begin{bmatrix}
	     \Sigma_x & 0 \\
	     0 & \Sigma_y
      \end{bmatrix}	  	  
	\right).
  \end{equation}
Let $\eye_m$ ($\eye_n$) be the identity matrix from $\mathbb{R}^{m \times m}$ ($\mathbb{R}^{n \times n}$). By Lemma \ref{thm:lttmnd},
\begin{equation}
  \begin{bmatrix}
    \eye_m & 0 \\ \mathbf{A} & \eye_n
  \end{bmatrix}\begin{pmatrix} \vec{x} \\ \vec{y} \end{pmatrix}  
  \label{eq:thmgaussianonestepproof}
\end{equation}
is a normal distribution with mean 
\begin{equation}
\begin{pmatrix}
  \mu_x \\ \mathbf{A} \mu_x
\end{pmatrix}
\end{equation}
and covariance matrix
\begin{equation}
  \begin{bmatrix}
    \eye_m & 0 \\ \mathbf{A} & \eye_n
  \end{bmatrix}
  \begin{bmatrix}
	     \Sigma_x & 0 \\
	     0 & \Sigma_y
  \end{bmatrix}
  \begin{bmatrix}
    \eye_m & \mathbf{A}' \\ 0 & \eye_n
  \end{bmatrix} = 	  \begin{bmatrix}
	     \Sigma_x & \Sigma_x \mathbf{A}' \\
	     \mathbf{A} \Sigma_x & \mathbf{A} \Sigma_x \mathbf{A}' + \Sigma_y
      \end{bmatrix},
\end{equation}
which finishes the proof.
\end{proof}

\subsection{Analysis of the EM-algorithm}

First, we will show that the hidden state and observation follow the normal distribution, and we will calculate its parameters.

\begin{lemma}
  \label{thm:normallydistributedattime}
  For an autonomous LDS $\LDS$ its hidden state follows the normal distribution
  \begin{equation}
    \hstatet \sim \normal\left(\hmatrix^t \hstatemu_0, \hmatrix^t \hstatecovz (\hmatrix')^t + \sum_{i = 0}^{t-1} \hmatrix^{i} \hcovariance (\hmatrix')^{i} \right),
    \label{eq:hstatedistro}
  \end{equation}
  and the observations follow the normal distribution
  \begin{equation}
    \obst \sim \normal\left(\omatrix\hmatrix^t \hstatemu_0, \omatrix\hmatrix^t \hstatecovz (\hmatrix')^t\omatrix' + \left[ \sum_{i = 0}^{t-1} \omatrix\hmatrix^{i} \hcovariance (\hmatrix')^{i}\omatrix' \right] + \ocovariance\right).
    \label{eq:obsdistro}
  \end{equation}
\end{lemma}
\begin{proof}
  For $t=0$, our assumption was that
  \begin{equation}
    \hstate_0 \sim \normal( \hstatemuz, \hstatecovz ),
    \label{eq:hstatezdistro}
  \end{equation}
  which proves \eqref{eq:hstatedistro} for $t=0$.
  
  The rest of the proof is done by the mathematical induction. Assume that $\hstatet$ follows the normal distribution stated in \eqref{eq:hstatedistro}. Then, according to Lemma \ref{thm:lttmnd}, $\hmatrix \hstatet$ follows normal distribution
  \begin{equation}
    \hmatrix \hstatet \sim \normal\left(\hmatrix \hmatrix^t \hstatemu_0, \hmatrix \left[ \hmatrix^t \hstatecovz (\hmatrix')^t + \sum_{i = 0}^{t-1} \hmatrix^{i} \hcovariance (\hmatrix')^{i} \right] \hmatrix' \right) =
    \normal\left(\hmatrix^{t+1} \hstatemu_0, \hmatrix^{t+1} \hstatecovz (\hmatrix')^{t+1} + \sum_{i = 1}^{t} \hmatrix^{i} \hcovariance (\hmatrix')^{i} \right).
  \end{equation}
  By Lemma \ref{thm:summnd},
  \begin{equation}
    \hstate_{t+1} = \hmatrix \hstatet + \hnoise_{t+1} \sim
    \normal\left(\hmatrix^{t+1} \hstatemu_0, \hmatrix^{t+1} \hstatecovz (\hmatrix')^{t+1} + \sum_{i = 0}^{t} \hmatrix^{i} \hcovariance (\hmatrix')^{i} \right),
  \end{equation}
  which finishes the proof. The proof for observation $\obst$ is analogous.
\end{proof}

As all the observations are normally distributed, we can ask whether their concatenation would be normally distributed as well. In that case, we might use algorithms for clustering a mixture of Gaussian distributions to cluster a mixture of LDS trajectories. We will answer this question in the next paragraphs.

\begin{lemma}
  \label{thm:concathiddennormal}
  Vector
  \begin{equation}
  \begin{pmatrix}
    \hstate_0 \\ \hstate_1 \\ \vdots \\ \hstate_\len
  \end{pmatrix}
  \end{equation}
  is normally distributed.
\end{lemma}
\begin{proof}
   The proof will be done by mathematical induction. Vector $\hstate_0$ is normally distributed by Assumption \ref{assumption:initialstate}.
   
   Assume that $(\hstate_0, \hstate_1, \ldots, \hstatet)'$ is normally distributed up to some time $\tii$. Then, as the noise is independent of the hidden state and between its samples (see Assumption \ref{assumption:independence}), $(\hstate_0, \hstate_1, \ldots, \hstatet, \hnoise_{t+1})'$ is normally distributed. By Lemma \ref{thm:gaussianonestep},
   \begin{equation}
      (\hstate_0, \hstate_1, \ldots, \hstatet, \hstate_{t+1})' = (\hstate_0, \hstate_1, \ldots, \hstatet, \hmatrix \hstatet + \hnoise_{t+1})'
   \end{equation}
   is normally distributed, where the transformation matrix $\mathbf{A}$ applied to vector $(\hstate_0, \hstate_1, \ldots, \hstatet, \hnoise_{t+1})'$ in Lemma \ref{thm:gaussianonestep} is equal to
   \begin{equation}
     \begin{bmatrix}
       \eye_n & 0 & \cdots & 0 & 0 \\
       0 & \eye_n & \cdots & 0 & 0 \\
       \vdots & \vdots & \ddots & \vdots & \vdots \\
       0 & 0 & \cdots & \hmatrix & \eye_n
     \end{bmatrix}.
   \end{equation}
   The proof is then finished by the standard mathematical induction argument.
%
\end{proof}
An alternative way to prove Lemma \ref{thm:concathiddennormal} would be to use a direct proof, similar to the proof of Lemma \ref{thm:gaussianonestep}. In that case, we can see that the transformation matrix needed to transform vector $(\hstate_0, \hnoise_1, \ldots, \hnoise_\len)'$ to $(\hstate_0, \hstate_1, \ldots, \hstate_\len)'$ is
\begin{equation}
   \label{eq:transformationmarix}
    \begin{bmatrix}
      \eye_\hdim & 0 & 0 & \cdots & 0 \\
      \hmatrix & \eye_\hdim & 0 & 0 & \cdots & 0 \\
      \hmatrix^2 & \hmatrix & \eye_\hdim & 0 & \cdots & 0 \\
      \hmatrix^3 & \hmatrix^2 & \hmatrix & \eye_\hdim & \cdots & 0 \\
      \vdots & \vdots & \vdots & \vdots & \ddots & \vdots \\
      \hmatrix^{\hdim-1} & \hmatrix^{\hdim-2} & \hmatrix^{\hdim-3} & \hmatrix^{\hdim-4} & \cdots & \eye_\hdim
    \end{bmatrix}.
\end{equation}
The linear transformation theorem \ref{thm:lttmnd} can then be used to calculate the exact parameters of the distribution.

\begin{corollary}
  \label{cor:concatobsnormal}
  Vector of concatenated observations
  \begin{equation}
    \begin{pmatrix}
      \obs_0 \\ \obs_1 \\ \vdots \\ \obs_\len
    \end{pmatrix}
  \end{equation}
  is normally distributed.
\end{corollary}
\begin{proof}
  By Lemma \ref{thm:concathiddennormal}, vector $(\hstate_0, \hstate_1, \ldots, \hstate_\len)'$ follows the normal distribution. By the linear transformation theorem $\omatrix (\hstate_0, \hstate_1, \ldots, \hstate_\len)'$ is normally distributed. Since the observation noise is independent of the state values and between its samples, Lemma \ref{thm:summnd} proves that 
  \begin{equation}
    \begin{pmatrix}
      \obs_0 \\ \obs_1 \\ \vdots \\ \obs_\len
    \end{pmatrix} =
    \omatrix'
    \begin{pmatrix}
      \hstate_0 \\ \hstate_1 \\ \vdots \\ \hstate_\len
    \end{pmatrix} + \begin{pmatrix}
      \onoise_0 \\ \onoise_1 \\ \vdots \\ \onoise_\len
    \end{pmatrix}
  \end{equation}
follows the normal distribution.
\end{proof}

Corollary \ref{cor:concatobsnormal} means that clustering a mixture of multiple LDSs is no more difficult than clustering a mixture of Gaussian distributions. We state this finding formally in the following theorem.
\begin{theorem}
   \label{thm:reduction1}
    There exists a polynomial reduction that reduces the problem of clustering a mixture of autonomous LDSs with hidden states to the clustering of a mixture of Gaussian distributions.
\end{theorem}
\begin{proof}
    The reduction comes from the Corollary \ref{cor:concatobsnormal}. In polynomial time, we can concatenate the vector of observations to a vector, one vector per trajectory. As the resulting concatenations are normally distributed, they can be clustered by any algorithm clustering a mixture of Gaussian distributions.
\end{proof}

Since there exists a reduction from clustering a mixture of autonomous LDS trajectories to clustering a mixture of Gaussian distributions, it is worth formally stating the reduction in the other way despite it being trivial to prove.

\begin{lemma}
\label{thm:reduction2}
    There exists a polynomial reduction from the problem of clustering a mixture of Gaussians to the clustering of a mixture of LDS trajectories.
\end{lemma}
\begin{proof}
    For a point in the Gaussian mixture, consider a trajectory with a length of $1$, where we set $\hdim = \odim$, $\omatrix = \eye_\hdim$, and let $\onoise = 0$ so that the observation is equal to the hidden state. For each point in the Gaussian mixture, we create a single trajectory of length $1$ where the initial hidden state $\hstate_0$ is set to equal the point. The problem of clustering the mixture of Gaussian distributions can then be solved by finding a clustering of a mixture of LDS trajectories, showing that the problem of clustering of LDS trajectories is at least as difficult as clustering a mixture of Gaussian distributions.
\end{proof}

\begin{theorem}
    Finding a soft clustering of a mixture of LDS trajectories with a log-likelihood within an additive factor of the optimal log-likelihood is NP-hard when $\clusters = 2$.
\end{theorem}
\begin{proof}
    The statement is a direct corollary of \ref{thm:reduction2}. The problem of clustering a mixture of Gaussian distributions is known to be NP-hard, even in the special case of spherical clusters. \cite{gaussiannphard} The initial conditions in proof of Lemma \ref{thm:reduction2} are defined so that the initial hidden state is propagated into the observation so that the original Gaussian distribution is clustered directly. Paper \cite{gaussiannphard} assumes that the variances are non-negligible and the Gaussians are spherical, which is a special case covered by the problem of clustering of LDSs. As the problem of clustering of LDSs includes a subset of inputs that can be used to solve an NP-hard problem, soft-clustering of LDSs is NP-hard.
\end{proof}

In the next section, we will focus on the consequences of the property that the concatenation of the observations is normally distributed. It is also worth mentioning that the result from Corollary \ref{cor:concatobsnormal} does not apply to LDSs with a control input as, in that case, the distribution cannot be modeled by only a single Gaussian, but a mixture of Gaussian distributions is needed (under similar assumptions). In the case of LDS with control input, Lemma \ref{thm:normallydistributedattime} does not hold.

\subsection{Implications of the Normally Distributed Observations}

As we have seen in the last section, finding the clustering of a mixture of autonomous LDSs is, in principle, the same as finding a clustering of a mixture of Gaussian distributions. As finding a clustering for a mixture of Gaussian distributions is a well-studied problem (and with more results than those that apply to the joint problem), we will summarize some of the important results in this section.

\begin{itemize}
  \item In general, the EM-algorithm is guaranteed to converge to a local minimum, maximum, or saddle point of the likelihood function under the assumption of continuity \cite{wuemalgo}. 
  \item The EM-algorithm is connected to gradient ascent. See paper \cite{jordanmixture} for details.
  \item If means of the Gaussians in the mixture are provided, local convergence to a global optimum of the likelihood function is guaranteed \cite{zhao2020}. The paper uses upper and lower bounds to prove that the EM algorithm update rule behaves as a contraction in the neighborhood of the global optimum.
  \item Paper \cite{NIPS2016_3875115b} shows that in the case of a mixture of more than two Gaussians, the local minima of the likelihood function can be arbitrarily bad, compared to the global optimum, even if the Gaussians are well-separated. The paper also gives a lower bound on convergence to bad critical points, which emphasizes the influence of the initialization on the final results.
  \item Recent paper \cite{pmlr-v125-kwon20a} proves a linear bound on the number of samples needed for EM-algorithm to converge in the case of a mixture of three or more spherical, well-separated Gaussians.
\end{itemize}


As can be seen, when there are three or more components in the mixture, the statistical guarantees are not favorable in the case of likelihood maximization using the EM-algorithm. Besides those general properties, when a mixture of only two Gaussians is considered, better convergence guarantees have been found in special cases.
\begin{itemize}
  \item Paper \cite{wu2021randomly} shows that with random initialization, the EM-algorithm form mixture of two Gaussians converges in $\mathcal{O}(\sqrt{n})$ with a high probability in Euclidean distance for sufficiently large $n$ (linearly growing with dimension up to a logarithmic factor). The result holds generally, even if no separation conditions are met.
  \item If we consider a mixture of two balanced Gaussians with known covariance matrices, there exist global convergence guarantees - given an infinite number of samples, the EM-algorithm converges geometrically to the correct mean vectors \cite{tensteps}.
  \item Paper \cite{NIPS2016_792c7b5a} proves convergence of the sequence of estimates for population EM when applied to a mixture of two Gaussians. The algorithm gives three possible optima for mean convergence and also provides parameter settings when the means are identified correctly or the algorithm converges to the point when the estimates are both the average of the true mean values. The results are then extended to the sample-based EM, and the probability of convergence is proven.
\end{itemize}
To contrast the previous paragraphs, even when there are two clusters with spherical Gaussians and shared variance, the soft clustering problem is NP-hard \cite{gaussiannphard}. The NP-hardness is proved for approximation of the log-likelihood within an additive factor. The same paper \cite{gaussiannphard} also shows that the NP-hardness remains for non-negligible variances. The complexity is shown by a reduction to the $k$-means problem.

Recent analyses focus on many special cases of the clustering of mixture of Gaussian distributions.
\begin{itemize}
    \item Paper \cite{dwivedi2020sharp} focuses on weakly identifiable models. They analyze the case of mixture of two equal-sized spherical Gaussian distributions that share the covariance matrices. The locations of the Gaussian distributions are then assumed to be symmetric with respect to axes origin. The paper than discusses the univariate case and shows that the statistical estimation error of the EM estimates is of the order of $n^{-\frac{1}{8}}$ and after $n^{\frac{3}{4}}$ steps, the error is in the order of $n^{-\frac{1}{4}}$. In the multivariate case, shared covariances improve the convergence criteria compared to the general case.
    \item Paper \cite{weinberger2022algorithm} studies a similar case - two symmetrically located spherical Gaussians, however, the mixture in this case is assumed imbalanced with known mixture coefficient. The authors then prove that the population-based EM-algorithm is globally convergent if the initial estimate has non-negative inner product with the mean of the larger component. When initialized to center the axis, error rate is given after a number of iterations inversely proportional to the mixing ratio and the norm of the cluster centers. Bounds for the empirical iteration are given as well.
    \item Further analyses of the weakly separated case are provided in \cite{ho2022weak}. The paper shows that the convergence rate is of the order of $n^{-\frac{1}{6}}$ or $n^{-\frac{1}{8}}$. The paper shows that sometimes the EM-algorithm shows high likelihood of the cluster means being equal despite this being false.
    \item In \cite{zhang2022distributed}, the authors develop a generalization of the standard EM-algorithm that can work in distributed setting. The method is consistent and retains the $\mathcal{O}(\sqrt{n})$ consistency under specified conditions. The authors then compare the method with some of the existing approaches, showing its superiority.
    \item Lastly, paper \cite{ho2022convergence} provides convergence rates for Gaussian mixtures of experts, which is a class of regression models. The authors state the notion of algebraic independence allowing them to establish a connection to partial differential equations, which in turn are used to prove the convergence rate. 
\end{itemize}

\subsection{Practical Applicability of the Gaussian Mixture-Based EM-algorithm}

Using EM-algorithm directly on concatenated vectors requires fitting $\Oh(\len (\odim + \hdim) \clusters)$ parameters in the case of mean values, and unfortunately, $\Oh(\len^2 (\hdim^2 + \odim^2) \clusters)$ parameters in the covariance matrix. By exploiting the transformation matrix in \eqref{eq:transformationmarix}, the number of parameters of the covariance matrices can be simplified by removing some degree of freedom from the problem, keeping only free parameters in $\hcovariance$, $\ocovariance$, $\hmatrix$, and $\omatrix$. Thus, we need only $\Oh((\hdim^2 + \odim^2) \clusters)$ parameters. {Adding those constraints can, however, cause loss of the theoretical properties of the EM-algorithm. 

A direct approach to solving the joint problem is to use the MLE estimates. In the case of spherical clusters with equal variance and under the negligence of the cost for the initial hidden state, the joint problem reduces to the minimization of 
\begin{equation}
  \min_{\hnoiset, \onoiset, \hstate_0, \hmatrix, \omatrix, \indicti} \sum_{i=1}^\trajectories \left( \sum_{\tii = 2}^\len \| \hnoiseti \|_2^2 + \sum_{\tii = 1}^\len \| \onoiseti \|_2^2 \right),
  \label{eq:lsqjoint}
\end{equation}
subject to
\begin{align}
  \hstateti &= (\hmatrix^{\indicti})\hstatetmoi + \hnoiseti, && \forall \tii \in \{2, 3, \ldots, \len\}, && \forall i \in \{1,2,\ldots,\trajectories\}, \\
  \obsti &= (\omatrix^{\indicti})'\hstateti + \onoiseti, && \forall \tii \in \{1, 2, \ldots, \len\}, && \forall i \in \{1,2,\ldots,\trajectories\}.
\end{align}
We can see in \eqref{eq:lsqjoint} that the MLE estimate requires to have a single parameter for each time step and each trajectory, which is the noise value assigned to the trajectory at a particular time. This means $\Oh(\len (m + n) \trajectories k)$ parameters, again too much for practical usability.

For completeness, the formula above leads to the following EM-heuristic formulation.

\begin{algorithm}
\begin{algorithmic}
\Function{EM-clustering}{$N$ trajectories $\obsti \in \mathbb{R}^{m \times T}$, $K$}
  \State\Comment{Generate a random partitioning into two clusters.}
  \State $\vec{l}_i \gets \Call{RandomInt}{\{0,1,\ldots, K - 1\}}$
  \State
  \State\Comment{Iterate until convergence.}
  \While{$l_i$ changes for any $i \in \{1,2, \ldots, N\}$}
     \State\Comment{For each cluster, find cluster parameters}
     \For{$c \in \{0,1,\ldots, K-1\}$}
        \State Find the cluster $C_{c}$ parameters by learning an LDS from a set of trajectories 
     \EndFor
    \State
     \State\Comment{Reassign the trajectories to the clusters.}
	\For{$i \in \{1, 2, \ldots N\}$}
       \State \begin{equation}\indicti \gets \displaystyle\mathop{\mathrm{arg\, min}}_{c \in \{0,1,\ldots, K-1\}}
		\min_{\hnoiseti, \onoiseti, \hstate_0} \left(      
       		\sum_{\tii = 2}^\len \| \hnoiseti \|_2^2 + \sum_{\tii = 1}^\len \| \onoiseti \|_2^2,  \right)  \label{eq:estep}
       		\end{equation} \\
       where each of the minimization problems is subject to
       \begin{align}
  \hstateti  &= (\hmatrix^{\ci})\hstatetmoi + \hnoiseti , && \forall \tii \in \{2, 3, \ldots, \len\}, \label{eq:estepcond1} \\
  \obsti &= (\omatrix^{\ci})'\hstateti + \onoiseti, && \forall \tii \in \{1, 2, \ldots, \len\},. \label{eq:estepcond2}
\end{align}
     \EndFor
  \EndWhile
\EndFunction
\end{algorithmic}
\caption{The EM heuristic with constraints.}
\label{alg:emheuristic}
\end{algorithm}
This algorithm is guaranteed to converge to a local optimum or a saddle point.

\subsection{Connection to $k$-means}

In our effort to improve the practical applicability of the algorithm, we can take inspiration from the mixture of the Gaussians approach. For the classical EM-algorithm, there exists a faster heuristic - Lloyd's algorithm \cite{lloydkmeans} for the $k$-means problem. In this section, we will show the connection of the minimization problem {from the main paper body} to the $k$-means problem and the connection of the heuristic to Lloyd's algorithm.

Recall the objective function,
\begin{equation}
    \min_{
      \substack{
       \vec{f}^0_t, \mathbf{F}_0, \mathbf{G}_0, \vec{\upsilon}^0_t, \vec{\omega}^0_t, \vec{\phi}_0^0, \\
       \vec{f}^1_t, \mathbf{F}_1, \mathbf{G}_1, \vec{\upsilon}^1_t, \vec{\omega}^1_t, \vec{\phi}_0^1 \\
       l_\tii
    }}
        \sum_{\tii=1}^N \sum_{t=1}^{T} 
        \|\vec{X}^\tii_t - \vec{f}^{l_i}_t\|_2^2 
    + 
       \sum_{c \in \{0, 1\}}\sum_{t=1}^{T} \left[  \| \vec{\upsilon_{t}}^{c} \|_2^2 + \| \vec{\omega_{t}}^{c} \|_2^2 \right].
    \label{eq:optcriterion}
\end{equation}
The first term of the objective function calculates the difference between the cluster means to the observations; the second term then minimizes noise that is induced by the optimal trajectory defined by the cluster means. With $\trajectories \to \infty$ the cost function goes to
\begin{equation}
    \min_{
      \substack{
       \vec{f}^0_t, \mathbf{F}_0, \mathbf{G}_0, \vec{\upsilon}^0_t, \vec{\omega}^0_t, \vec{\phi}_0^0, \\
       \vec{f}^1_t, \mathbf{F}_1, \mathbf{G}_1, \vec{\upsilon}^1_t, \vec{\omega}^1_t, \vec{\phi}_0^1 \\
       l_\tii
    }}
        \sum_{\tii=1}^N \sum_{t=1}^{T} 
        \|\vec{X}^\tii_t - \vec{f}^{l_i}_t\|_2^2 
    ,
    \label{eq:optcriterion2}
\end{equation}
as the other terms do not increase with the number of trajectories. The formula in \eqref{eq:optcriterion2} is the standard $k$-means criterion. Applying the same reasoning to the EM-heuristic {in the main paper body}, leads to the standard Lloyd's algorithm, as with $\trajectories \to \infty$, minimization
\begin{equation}
    \min_{
      \substack{
       \vec{f}^c_t, \mat{F}_c, \mat{G}_c, \vec{\upsilon}^c_t, \vec{\omega}^c_t, \vec{\phi}_c^0
    }}
      \left[ \sum_{\tii=1}^N \sum_{t=1}^{T} 
        \mathbbm{1}[l_i = c] \cdot
        \|\vec{X}^\tii_t - \vec{f}^{c}_t\|_2^2 \right] +
        \| \vec{\upsilon}^{c} \|_2^2 + \| \vec{\omega
}^{c} \|_2^2 
\end{equation}
converges to the following minimization problem
\begin{equation}
  \min_{
      \substack{
       \vec{f}^c_t, \mat{F}_c, \mat{G}_c, \vec{\upsilon}^c_t, \vec{\omega}^c_t, \vec{\phi}_c^0
    }}
      \left[ \sum_{\tii=1}^N \sum_{t=1}^{T} 
        \mathbbm{1}[l_i = c] \cdot
        \|\vec{X}^\tii_t - \vec{f}^{c
}_t\|_2^2 \right],
\end{equation}
which is minimized by $f_t^c$ being the cluster means. To wrap this up, with the increasing number of trajectories $\trajectories \to \infty$, the EM-heuristic converges to Lloyd's algorithm \cite{lloydkmeans} for the $k$-means problem.

\section{Clustering Performance Metrics}
\label{sec:metrics}
The $F_1$ score, which has been widely used in classification performance measurements, is defined as 
\begin{align}
    F_1 =2 \cdot \frac{\textrm{precision} \cdot \textrm{recall}}{\textrm{precision}+ \textrm{recall}}
\end{align}
where precision $=\frac{TP}{TP+FP}$ and recall $=\frac{TP}{TP+FN}$ and $TP$, $FP$, and $FN$ are the numbers of true positives, false positives, and false negatives, respectively.
We calculate the $F_1$ score twice for each class, once with one class labeled as positive and once with the other class labeled as positive, and then we select the higher score for each class. This approach is used because there is no predefined positive and negative labels.

\end{document}